\documentclass[12pt]{amsart}
\usepackage{amsmath}
\usepackage{amsthm}
\usepackage{amssymb}
\usepackage{amsfonts,mathrsfs}
\usepackage{stmaryrd}
\usepackage{amsxtra}  
\usepackage{epsfig}
\usepackage{verbatim}

\theoremstyle{plain}
\newtheorem{theorem}[equation]{Theorem}
\newtheorem{proposition}[equation]{Proposition}
\newtheorem{lemma}[equation]{Lemma}

\newtheorem{definition}[equation]{Definition}
\theoremstyle{remark}

\numberwithin{equation}{section}

\include{header}

\newcommand{\dbarstar}{\bar{\partial}^{\star}}
\newcommand{\ssubset}{\subset\subset}
\newcommand{\psum}{\sideset{}{^{\prime}}{\sum}}
\newcommand{\sn}{{\mathscr N}}
\newcommand{\dbar}{\bar \partial}

\begin{document}

\bibliographystyle{plain}
\title[Regularity of the Bergman projection]{Regularity of the Bergman projection on forms
and plurisubharmonicity conditions}
\author{A.-K. Herbig \& J.D. McNeal}
\subjclass{32W05}
\thanks{Research of the first author was partially supported by a Rackham Fellowship}
\thanks{Research  of the second author was partially supported by an NSF grant}
\address{Department of Mathematics, \newline University of Michigan, Ann Arbor, Michigan 48109}
\email{herbig@umich.edu}
\address{Department of Mathematics, \newline Ohio State University, Columbus, Ohio 43210}
\email{mcneal@math.ohio-state.edu}
\date{}
\begin{abstract}
  Let $\Omega\ssubset\mathbb{C}^{n}$ be a smoothly bounded domain. Suppose $\Omega$ has a  
  defining function, such that the sum of any $q$ eigenvalues of its complex Hessian is 
  non-negative. We show that this implies global regularity of the Bergman projection, 
  $B_{j-1}$, and the $\dbar$-Neumann operator, $N_{j}$, acting on $(0,j)$-forms, for $j\in\{q,\dots,n\}$.
\end{abstract}
\keywords{Bergman projection, global regularity}
\maketitle
\section{Introduction}
A function $f\in C^{\infty}(\Omega)$ is holomorphic on $\Omega$, if it satisfies the Cauchy-Riemann equations: $\dbar f=\sum_{k=1}^{n}\frac{\partial f}{\partial\bar{z}_{k}}\, d\bar z_k=0$ in $\Omega$.
Denote the set of holomorphic functions on $\Omega$ by $H(\Omega)$.
The Bergman projection, $B_{0}$, is the orthogonal projection of square-integrable functions onto $H(\Omega)\cap L^{2}(\Omega)$. Since the Cauchy-Riemann operator, $\dbar$ above, extends naturally to act on higher order forms, we can as well define Bergman projections on higher order forms:  let $B_j$ be the orthogonal projection of square-integrable $(0,j)$-forms onto its subspace of $\dbar$-closed, square-integrable $(0,j)$-forms. 

In this paper we give a condition on $\Omega$ which implies that $B_j f$ is smooth on $\overline{\Omega}$ whenever f is. Our result is the following:
\begin{theorem}\label{T:MainT}
  Let $\Omega\ssubset\mathbb{C}^{n}$ be a smoothly bounded domain. Suppose there exists a smooth  
  defining function of $\Omega$, such that the sum of any $q$ eigenvalues of its complex Hessian is 
  non-negative on $\overline\Omega$. Then the Bergman projection, $B_{j-1}$, is globally regular for $q\leq j\leq n$.
\end{theorem}
Global regularity of the Bergman projection is closely tied to the regularity of the $\dbar$-Neumann operator. Recall that the $\dbar$-Neumann operator, $N_{j}$, is the operator, acting on square-integrable $(0,j)$-forms, which inverts a particular boundary value problem associated to the complex Laplacian. 
By a result of Boas and Straube \cite{Boa-Str90}, $N_{j}$ is globally regular if and only if $B_{j-1}$, $B_{j}$ and $B_{j+1}$ are. Consequently, the hypothesis of our Theorem \ref{T:MainT} implies that $N_{j}$, $q\leq j\leq n$, are also globally regular.

Theorem  \ref{T:MainT} is an extension of an earlier result of Boas and Straube in
\cite{Boa-Str91}, though only a partial one.
There the authors show that $B_j$ is globally regular for all $j\in\{1,\dots,n\}$, if $\Omega$ has a smooth defining function which is plurisubharmonic on $b\Omega$ (the boundary of the domain $\Omega$). Our theorem covers the case considered by Boas-Straube, when $q=1$, but only under the stronger hypothesis that the defining function $r$ is plurisubharmonic on all $\overline\Omega$.
Our method of proof is quite different than the one in \cite{Boa-Str91}, though there are, naturally, some points in  common. Our proof shares more similarities with one of Kohn \cite{Koh99}, where he determined how the range of Sobolev norms $\|\cdot\|_k$, where $\| B_0f\|_k\leq C\| f\|_k$ holds, depends on the Diederich-Forn\ae ss exponent -- the (largest) exponent $1>s>0$ such that $-(-r)^s$ is plurisubharmonic -- when $r$ itself is not plurisubharmonic.

To compare our proof with that of \cite{Boa-Str91}, consider the case in common to both results, i.e., assume that $r$ is plurisubharmonic on all of $\overline\Omega$. The essential problem is to estimate some Sobolev norm larger than $1/2$ of $B_*f$ by the same Sobolev norm of $f$, say $\| B_* f\|_1$ by
$\| f\|_1$. In both proofs, standard arguments reduce this problem to that of estimating $\| XB_* f\|_0$ by
$\| f\|_1$, where $X$ is a tangential vector field to $b\Omega$  which is tranverse to the complex tangent space. In order to achieve this estimate, $X$ must commute ``nicely'' with the $\dbar$-complex in some fashion, so that one can absorb the error terms which arise in comparing 
$X\left(B_*f\right)$ to $B_*\left(Xf\right)$. In \cite{Boa-Str91}, the focus is on the commutator $\left[X, \dbar\right]$, and the plurisubharmonicity of $r$ is used to construct a special vector field $X$ so that this commutator has a small component in the 
complex normal direction to $b\Omega$. The role of plurisubharmonicity in this approach is that non-negativity of the matrix $\left(\frac{\partial^2 r}{\partial z_k\partial\bar{z}_l}\right)$ 
can be used to adjust any tangential, tranverse field by adding tangential fields containing only barred derivatives (which are benign in the estimates considered) to it in such a way that the commutator has the desired property. In our proof, the focus is on the commutator $\left[\dbar^*, X\right]$, for $X$ the ``natural'' tangential vector field tranverse to the complex tangent space (see the subsection I Section 4 below). Actually, we focus on a tangential field $T$, very closely related to $X$; the crucial property
of $T$ is that it preserves the domain of $\dbarstar$. We point out that the passage from $X$ to $T$
is a lower-order adjustment and does not depend on the plurisubharmonicity of $r$, i.e., $T$ differs from
$X$ by a 0th order operator.
The matrix $\left(\frac{\partial^2 r}{\partial z_k\partial\bar{z}_l}\right)$ then appears as the matrix of coefficients in front of $T$, acting on various components of a $(0,q)$-form $\varphi$, closely connected to the form $f$, paired with a neglible form. 
The non-negativity of $\left(\frac{\partial^2 r}{\partial z_k\partial\bar{z}_l}\right)$ then allows the use of the Cauchy-Schwarz inequality to separate this pairing into separate, purely quadratic factors, see \eqref{E:CauSch} below. The factor involving the $T$ derivatives of $\varphi$ is then estimated by a small constant times the $\dbar$-Dirichlet form of $\varphi$ using the $\dbar$-Hardy inequality proved in Section \ref{S:BasEst}.

Theorem \ref{T:MainT} gives many examples of domains where the Bergman projection on higher-level forms is regular while the Bergman projection on functions is not.
Suppose $D\ssubset\mathbb{C}^{n}$ is a smoothly bounded domain, and $\rho$ a smooth defining function of $D$.  Let $C$ be a lower bound for the sum of the eigenvalues of $i\partial\dbar\rho$ on $\overline{D}$, and define
\begin{align*}
  \widetilde{D}=\{(z,w)\in\mathbb{C}^{n+m}\;|\;
  r(z,w):=\rho(z)+K(|w_{1}|^2+\dots+|w_{m}|^{2})<0\},
\end{align*}
where $K\geq 0$ is chosen such that $K\geq|C|$. Then $\widetilde{D}$ is a smoothly bounded domain, and the sum of any $(n+1)$ eigenvalues of $i\partial\dbar r$ is non-negative on the closure of $\widetilde{D}$. Thus, by Theorem \ref{T:MainT}, the Bergman projection, 
$B_{q}^{\widetilde{D}}$, on $(0,q)$-forms is globally regular for $n\leq q\leq n+m$. 
In \cite{Bar84}, Barrett constructed a smoothly bounded domain $D'$ in $\mathbb{C}^{2}$ for which the Bergman projection on functions fails to be regular. Inserting Barrett's domain $D'$ for $D$ 
in the preceding construction, one obtains a smoothly bounded domain $\widetilde{D'}$ in $\mathbb{C}^{2+m}$ such that $B_{q}^{\widetilde{D'}}$ is regular for $2\leq q\leq 2+m$. However, using similar arguments to those in \cite{Bar84}, one can show that $B_{0}^{\widetilde{D'}}$ fails to be regular.

The paper is structured as follows. In Section \ref{S:Prelim} we present the general setting and a brief review of the $\dbar$-Neumann problem and its  relation to the Bergman projections. In Section \ref{S:BasEst} we derive the basic estimates which will be used for the proof of Theorem \ref{T:MainT}. In particular, we prove a Hardy-like inequality for the $\dbar$-complex, Proposition \ref{P:dbarHardy};
this estimate is of independent interest and should have further application.
In Section \ref{S:TheProof},  
we give the proof of Theorem \ref{T:MainT}. Although the total length of this paper exceeds that
of \cite{Boa-Str91},
the analytic heart of our proof is relatively short and is labeled as such in Section \ref{S:TheProof}.

We would like to thank K. Koenig  and E. Straube for pointing out  errors in an earlier version, and for their helpful suggestions.

\section{Preliminaries}\label{S:Prelim}
Throughout, let $\Omega\ssubset\mathbb{C}^{n}$ be a smoothly bounded domain, i.e., $\Omega$ is bounded and there exists a smooth, real-valued function $r$ such that $\Omega=\{z\in\mathbb{C}^{n}\;|\;r(z)<0\}$, and $\nabla r\neq 0$ when $r=0$. The hypothesis of Theorem \ref{T:MainT} will be abbreviated as follows
\begin{definition}\label{Hq} We say that a smoothly bounded domain $\Omega\ssubset\mathbb{C}^n$ satisfies
condition $(H_q)$ if there exists a defining function $r$ for $\Omega$ such that the sum of any $q$ eigenvalues of $\left(\frac{\partial^2 r}{\partial z_k\partial\bar{z}_l}\right)$ is non-negative on $\overline\Omega$.
\end{definition}

We shall write an arbitrary $(0,q)$-form $u$,  $0\leq q\leq n$, as
\begin{align}\label{E:Defu}
  u=\psum_{|J|=q}u_{J}d\bar{z}^{J},
\end{align}
where \begin{tiny}$\psum_{|J|=q}$\end{tiny} means that the sum is taken over strictly increasing multi-indices $J$ of length $q$. We define the coefficients $u_I$ for arbitrary multi-indices $I$ of length $q$, so that the $u_I$'s are antisymmetric functions in $I$.

Let $\Lambda^{0,q}(\overline{\Omega})$ and $\Lambda_{c}^{0,q}(\Omega)$ denote the $(0,q)$-forms with coefficients in 
$C^{\infty}(\overline{\Omega})$ and $C_{c}^{\infty}(\Omega)$, respectively. For $(0,1)$-forms, we use the pointwise inner product $\langle.,.\rangle$ defined by $\langle d\bar{z}^{k},d\bar{z}^{l}\rangle=\delta_{l}^{k}$. By linearity we extend this inner product to $(0,q)$-forms. We define the global $L^{2}$-inner product on $\Omega$ by
\begin{align*}
  (u,v)_{\Omega}=\int_{\Omega}\langle u,v\rangle\;dV\;\;\text{for}\;\;u,v\in\Lambda_{c}^{0,q}(\Omega),
\end{align*}
where $dV$ is the euclidean volume form. The $L^{2}$-norm of $u\in\Lambda_{c}^{0,q}(\Omega)$ is then given by $\|u\|_{\Omega}=\sqrt{(u,u)}_{\Omega}$, and we define the space $L_{0,q}^{2}(\Omega)$ to be the completion of $\Lambda_{c}^{0,q}(\Omega)$ under $\|.\|_{\Omega}$.

The Cauchy-Riemann operator, $\dbar$, acting on $u\in\Lambda^{0,q}(\overline{\Omega})$ is defined as follows
\begin{align*}
  \dbar u=\psum_{|J|=q}\sum_{k=1}^{n}\frac{\partial u_{J}}{\partial\bar{z}_{k}}d\bar{z}^{k}\wedge 
  d\bar{z}^{J},
\end{align*} 
where $u$ is expressed as in (\ref{E:Defu}). In order to be able to use Hilbert space techniques, we want to extend $\dbar$ to act on (a dense subspace of) $L^{2}_{0,q}(\Omega)$. To do this, first
extend $\dbar$ to act on non-smooth forms in the sense of distributions. Then, to obtain
a Hilbert space operator, restrict the domain of $\dbar$ to those forms $g\in L^{2}_{0,q}(\Omega)$, such that $\dbar g$, in the sense of distributions, is in $L_{0,q+1}^{2}(\Omega)$. In this way, $\dbar$ is a closed, densely defined operator on Hilbert spaces. We define its $L^{2}$-adjoint $\dbarstar$ in the usual manner:
a form $v\in L^{2}_{0,q+1}(\Omega)$ belongs to the domain of $\dbarstar$, $\text{Dom}(\dbarstar)$, if there exists a constant $C>0$ so that
\begin{align}\label{D:DomdbarstarDef}
  |(\dbar u, v)|\leq C\|u\|\;\;\;\text{holds for all} \;\;u\in\text{Dom}(\dbar).
\end{align}
The Riesz representation theorem implies that, if $v\in\text{Dom}(\dbarstar)$, there exists a unique $w\in L_{0,q}^{2}(\Omega)$, such that
\begin{align*}
  (\dbar u,v)=(u,w)
\end{align*}
holds for all $u\in\text{Dom}(\dbar)$; we write $\dbarstar v$ for $w$. One can show, using integration by parts, that if $v\in\mathcal{D}^{0,q}(\Omega):=\text{Dom}(\dbarstar)\cap\Lambda^{0,q}(\overline{\Omega})$, then $v$ satisfies the following boundary conditions:
\begin{align*}
  \sum_{k=1}^{n}\frac{\partial r}{\partial z_{k}}v_{kI}=0\;\;\text{on}\;\;b\Omega
\end{align*}
for any strictly increasing multi-index $I$ of length $q-1$. Here we mean by $kI$ the multi-index $\{k,I\}$. Denote by $\text{Dom}(\Box_{q})$ those $(0,q)$-forms $u\in\text{Dom}(\dbarstar)\cap\text{Dom}(\dbar)$ for which $\dbar u\in\text{Dom}(\dbarstar)$ and $\dbarstar u\in\text{Dom}(\dbar)$ holds. The operator $\Box_{q}=\dbar\dbarstar+\dbarstar\dbar$, defined for forms in $\text{Dom}(\Box_{q})$, is called the complex Laplacian.

We introduce a convenient, if non-standard, piece of notation: if $f$ is a $C^2$ function
\begin{align}\label{hess_q}
i\partial\dbar f(u,u):=\psum_{|J|=q-1}\sum_{k,l=1}^{n}\frac{\partial^2 f}{\partial z_{l}\partial\bar{z}_{k}}u_{lJ}
\bar{u}_{kJ},\qquad u\in\Lambda^{0,q}(\overline{\Omega}).
\end{align}
When $q=1$, \eqref{hess_q} is standard notation and expresses the natural action of the $(1,1)$-form $i\partial\dbar f$ on the pair of vectors $u$ and $\bar u$. For $q>1$, the left-hand side of \eqref{hess_q} does not have such a natural meaning. However, the right-hand side of \eqref{hess_q} arises repeatedly in integration by parts arguments on the $\dbar$ complex, and it is useful to abbreviate this expression by the left-hand side of \eqref{hess_q} for all levels of forms. For example, the basic identity for the $\dbar$-Neumann problem assumes the following form: if $u\in\mathcal{D}^{0,q}(\Omega)$,

\begin{align}\label{E:BasId}
\left\|\dbar u\right\|^2+\left\|\dbarstar u\right\|^2 =
\psum_{|I|=q}\sum_{k=1}^{n}\left\|\frac{\partial u_{I}}{\partial\bar{z}_{k}}\right\|^2+\int_{b\Omega}
i\partial\dbar r(u,u).
\end{align}
We also mention the equivalence of the following two facts related to \eqref{hess_q}:
\begin{enumerate}
\item[(i)] $i\partial\dbar f(u,u)\geq C|u|^2$ for all $ u\in\Lambda^{0,q}$.
\item[(ii)] The sum of any $q$ eigenvalues of the matrix $\left(\frac{\partial^2 f}{\partial z_k\partial \bar{z}_l}\right)$ is greater than or equal to $C$.
\end{enumerate}
A proof of the equivalence of (i) and (ii) follows by diagonalizing the matrix 
$\left(\frac{\partial^2 f}{\partial z_k\partial\bar{z}_l}\right)$; see \cite{Hor65} or \cite{Cat86}.

Suppose that for all  $u\in\mathcal{D}^{0,q}(\Omega)$, 
$i\partial\dbar r(u,u)\geq 0$ on the boundary of $\Omega$. Starting with \eqref{E:BasId}, one can show that
\begin{align}\label{E:VeryBasEst}
\|u\|^2+\psum_{|I|=q}\sum_{k=1}^{n}\left\|\frac{\partial u_{I}}{\partial\bar{z}_{k}}\right\|^2
\leq C\left(\left\|\dbar u\right\|^2+\left\|\dbarstar u\right\|^2\right)
\end{align}
holds for all $u\in\mathcal{D}^{0,q}(\Omega)$; here $C>0$ does not depend on $u$. If inequality 
(\ref{E:VeryBasEst}) holds, then the $\dbar$-Neumann operator exists: $N_{q}:L_{0,q}^{2}(\Omega)\longrightarrow\text{Dom}(\Box_{q})$ with $N_{q}\Box_{q}=Id$ on $\text{Dom}(\Box_{q})$ and $\Box_{q} N_{q}=Id$. One of the equations which connects the $\dbar$-Neumann operator and the Bergman projection is Kohn's formula, which says
\begin{align*}
  B_{q-1}=Id-\dbarstar N_{q}\dbar.
\end{align*}

Throughout the paper, we shall use the notation $|A|\lesssim|B|$ to mean
$|A|\leq C|B|$ for some constant $C>0$, which is independent of relevant parameters. It will be mentioned, or clear from the context, what those parameters are. We call the often-used inequality $|AB|\leq\eta A^{2}+\frac{1}{4\eta}B^{2}$ for $\eta>0$ the (sc)-(lc) inequality. Finally, we denote the commutator of two operators, $L$ and $M$, as usual: $\left[ L,M\right] =LM -ML$.

\section{Basic Estimates}\label{S:BasEst}
In this section, we derive a Hardy-like inequality for the $\dbar$-complex. This inequality, (\ref{E:dbarHardy}), is essential for our proof of Theorem \ref{T:MainT}. We start out with an energy identity for the $\dbar$-complex.
\begin{proposition}\label{P:BasEst}
Let $\Omega\ssubset\mathbb{C}^{n}$ be a smoothly bounded domain, $r$ a smooth defining function of $\Omega$. Let $s\geq 0$ and set $\tau=(-r)^s$. Then
\begin{align}\label{E:BasEst}
   \left \|\sqrt{\tau}\dbar u\right\|^{2}&+\left\|\sqrt{\tau}\dbarstar u\right\|^{2}-\int_{b\Omega}\tau i\partial\dbar r(u,u)\\
    & =
    \psum_{|I|=q}\sum_{k=1}^{n}\left\|\sqrt{\tau}\frac{\partial u_I}{\partial\bar{z}_{k}}\right\|^{2}
    -\int_{\Omega}i\partial\dbar\tau(u,u)
   +2\text{Re}\left([\dbarstar,\tau]u,\dbarstar u\right)\notag
\end{align} 
holds for all $u\in\mathcal{D}^{0,q}(\Omega)$.
\end{proposition}
Equation (\ref{E:BasEst}) was proved in \cite{McN02} for any $\tau\in C^2(\overline\Omega)$.
Other, related identities, e.g., for $(0,1)$-forms and forms vanishing on $b\Omega$, have been
obtained by several authors, starting with the basic work of Ohsawa and Takegoshi \cite{OhsTak87};
see \cite{McN05} for references and an expository account of these identities.
 However, since we need the identity when $\tau=(-r)^s$, for ranges of $s$ for which $\tau\notin  C^2(\overline\Omega)$, 
we give the modification of the proof in \cite{McN02} which yields Proposition \ref{P:BasEst}.
\begin{proof}
  Let $J,\;M$ and $N$ be multi-indices with $|J|=q-1$ and $|M|=|N|=q$. For notational ease denote 
  $\frac{\partial}{\partial z_{k}}$ by $\partial_k$. We write for $u\in\mathcal{D}^{0,q}(\Omega)$
  \begin{align*}
    \dbar u=\psum_M\sum_{k=1}^{n}\dbar_k u_M d\bar{z}^{k}\wedge d\bar{z}^M\;\;\text{and}\;\;
    \dbarstar u=-\psum_J\sum_{l=1}^{n}\partial_l u_{lJ}d\bar{z}^{J}.
  \end{align*}
  Then we obtain
  \begin{align*}
    \left\|\sqrt{\tau}\dbar u\right\|^2+\left\|\sqrt{\tau}\dbarstar u\right\|^2
    =&\psum_{M,N}\sum_{k,l=1}^{n}\sigma_{lN}^{kM}\left(\tau\dbar_{k} u_{M},\dbar_{l}u_{N}\right)\\
    &+\psum_{J}\sum_{k,l=1}^{n}\left(\tau\partial_{k}u_{kJ},\partial_{l}u_{lJ}\right),
     \end{align*}
     where $\sigma_{lN}^{kM}$  is the sign of the permutation 
     $\left(
     \begin{smallmatrix}
     kM\\ lN
     \end{smallmatrix}
     \right)$
     and equals zero whenever $\{k\}\cup\{M\}\neq \{l\}\cup\{N\}$.  Rearranging terms we obtain
     \begin{align*}
       \left\|\sqrt{\tau}\dbar u\right\|^2+\left\|\sqrt{\tau}\dbarstar u\right\|^2
       =&
       \psum_{M}\sum_{k=1}^{n}\left\|\sqrt{\tau}\dbar_{k} u_{M}\right\|^2\\
       &+\psum_{J}\sum_{k,l=1}^{n}\int_{\Omega}\tau\left(\partial_{l}u_{lJ}\bar{\partial_{k}}\bar{u}_{kJ}
       -\dbar_{k}u_{lJ}\partial_{l}\bar{u}_{kJ}\right),
     \end{align*}
     where we denote the last term on the right hand side by (I). Integration by parts gives
     \begin{align*}
       (I)=\psum_{J}\sum_{k,l=1}\int_{\Omega}\partial_{l}\left(\tau\dbar_{k} u_{lJ}\right)\bar{u}_{kJ}
       -\dbar_{k}\left(\tau\partial_{l}u_{lJ}\right)\bar{u}_{kJ},
     \end{align*}
     here no boundary integral appears because $\tau=(-r)^s$ is zero on the boundary of $\Omega$. 
     Since $\partial_{l}\dbar_{k}=\dbar_{k}\partial_{l}$ it follows
     \begin{align*}
     (I)=\psum_{J}\sum_{k,l=1}^{n}\int_{\Omega}\partial_{l}(\tau)\left(\dbar_{k}u_{lJ}\right)\bar{u}_{kJ}
     -\dbar_{k}(\tau)\left(\partial_{l}u_{lJ}\right)\bar{u}_{kJ}.
     \end{align*} 
     We would like to integrate the first term on the right hand side of the above equation by parts, but some care has to be taken since $\partial_{l}\tau$ is not defined on $b\Omega$ for $s\in(0,1)$. For $\epsilon>0$ small set $\Omega_{\epsilon}=\left\{z\in\Omega\,|-\epsilon<r(z)<0\right\}$. Then
\begin{align*}
  \psum_{J}\sum_{k,l=1}^{n}\int_{\Omega_{\epsilon}}\partial_{l}(\tau)\left(\dbar_{k}u_{lJ}\right)\bar{u}_{kJ}
  =&-\psum_{J}\sum_{k,l=1}^{n}\int_{\Omega_{\epsilon}}\dbar_{k}\bigl(\partial_{l}(\tau)\bar{u}_{kJ}\bigr)u_{lJ}\\
  &+\psum_{J}\sum_{k,l=1}^{n}\int_{b\Omega_{\epsilon}}\partial_{l}(\tau)u_{lJ}\dbar_{k}(r)\bar{u}_{kJ}
  \frac{dS}{|\partial r|}.
\end{align*}
Let us denote the boundary integral by (II). Using that $\tau=(-r)^s$ we can express (II) in the following manner
\begin{align*}
  (II)&=\psum_{J}\sum_{k,l=1}^{n}\int_{b\Omega_{\epsilon}}-s(-r)^{s-1}\partial_{l}(r)u_{lJ}\dbar_{k}(r)
  \bar{u}_{kJ}\frac{dS}{|\partial r|}\\
  &=-\psum_{J}\int_{b\Omega_{\epsilon}}s\epsilon^{s-1}\biggl|\sum_{l=1}^{n}
  \partial_{l}(r) u_{lJ}\biggr|^2\frac{dS}{|\partial r|}.
\end{align*}
Recall that $u\in\mathcal{D}^{0,q}(\Omega)$ means that $\sum_{l=1}^{n}\partial_{l}(r)u_{lJ}=0$ on $b\Omega$ for all increasing multi-indices $J$. Thus $\sum_{l=1}^{n}\partial_{l}(r)u_{lJ}
=O(\epsilon)$ on $b\Omega_{\epsilon}$, which yields $(II)=O(\epsilon^{s+1})$. Therefore, taking the limit as $\epsilon$ approaches $0$, it follows that
\begin{align*}
  \psum_{J}\sum_{k,l=1}^{n}\int_{\Omega}\partial_{l}(\tau)\left(\dbar_{k}u_{lJ}\right)\bar{u}_{kJ}
  &=-\psum_{J}\sum_{k,l=1}^{n}\int_{\Omega}\dbar_{k}\bigl(\partial_{l}(\tau)\bar{u}_{kJ}\bigr)u_{lJ}\\
  &=-\int_{\Omega}i\partial\dbar\tau (u,u)-\psum_{J}\sum_{k,l=1}^{n}\int_{\Omega}\partial_{l}(\tau)u_{lJ}
  \dbar_{k}\bar{u}_{kJ},
\end{align*}
which proves our claimed equation (\ref{E:BasEst}).
\end{proof}
In order to prove our main result of this section, inequality (\ref{E:dbarHardy}),  we need to show that under the hypotheses of Theorem \ref{T:MainT}, one can choose a defining function, $r$, for $\Omega$ with a suitable lower bound on the complex Hessian of $-(-r)^s$, $s\in(0,1)$. 
\begin{lemma}\label{L:wlog}
If $\Omega$ satisfies condition $(H_q)$, then for each $s\in(0,1)$, there exists a smooth defining function $r$ satisfying
  \begin{equation*}\label{E:wlog}
    i\partial\dbar\bigl(-(-r)^{s}\bigr)(u,u)\gtrsim(-r)^{s}|u|^2\;\;\text{for all}\;\;u\in\Lambda^{0,q}(\overline{\Omega}).
  \end{equation*}
  The constant in $\gtrsim$ depends on $s$, but is independent of $z\in\Omega$ and
  $u\in\Lambda^{0,q}$.
\end{lemma}

For $q=1$, Lemma \ref{L:wlog} was proved in \cite{Die-For77} (also see \cite{Ran81}), though it
was not stated in this form. The proof for general $q$ follows the same lines.

\begin{proof} Let $\rho$ satisfy Definition \ref{hess_q} and set $r(z):=e^{-K|z|^2}\rho(z)$, for 
a constant $K>0$ to be determined.

  A straightforward computation gives
  \begin{align*}
    i\partial\dbar\bigl(-(-r)^{s}\bigr)(u,u)=&s(-\rho)^{s-2}e^{-sK|z|^2}
    \left\{K\rho^{2}\Biggl[q|u|^{2}-sK\psum_{|J|=q-1}\biggl|\sum_{k=1}^{n}\bar{z}_{k}u_{kJ}\biggr|^{2}\Biggr]\right.\\
    &-\rho\Biggl[i\partial\dbar\rho(u,u)-2sK\text{Re\;}\biggl(\psum_{|J|=q-1}\sum_{k,l=1}^{n}z_{k}\bar{u}_{kJ}
    \frac{\partial\rho}{\partial z_{l}}u_{lJ}\biggr)\Biggr]\\
    &\left.+(1-s)\psum_{|J|=q-1}\biggl|\sum_{k=1}^{n}\frac{\partial \rho}{\partial z_{k}}u_{kJ}\biggr|^{2} 
    \right\}.
  \end{align*}
  Note that by (sc)-(lc) inequality we have
  \begin{align*}
    2sK\rho\text{Re}\;\biggl(\sum_{k,l=1}^{n}z_{k}\bar{u}_{kJ}
    \frac{\partial\rho}{\partial z_{l}}u_{lJ}\biggr)
    \geq
    (s-1)\biggl|\sum_{k=1}^{n}\frac{\partial \rho}{\partial z_{k}}u_{kJ}\biggr|^{2}
    -\frac{(sK\rho)^{2}}{1-s}\biggl|\sum_{k=1}^{n}\bar{z}_{k}u_{kJ}\biggr|^{2}
  \end{align*}
  for all increasing multi-indices $J$ with $|J|=q-1$.
  Using this inequality and that $i\partial\dbar\rho\geq 0$, it follows
  \begin{align*}
    i\partial\dbar\bigl(-(-r)^s\bigr)(u,u)
    \geq
    sK(-r)^s\biggl\{q|u|^{2}-K\frac{s}{1-s}\psum_{|J|=q-1}\biggl|\sum_{k=1}^{n}\bar{z}_{k}u_{kJ}\biggr|^{2}\biggr\}. 
  \end{align*}   
  Since $\Omega$ is a bounded domain, there exists a constant $D>0$ such that $|z|^{2}\leq D$ for  
  $z\in\overline{\Omega}$. We obtain
  \begin{align*}
    \psum_{|J|=q-1}\biggl|\sum_{k=1}^{n}\bar{z}_{k}u_{kJ}\biggr|^{2}\leq n q |z|^{2}|u|^{2}\leq n q D|u|^{2},
  \end{align*}
  which implies that
  \begin{align*}
    i\partial\dbar\bigl(-(-r)^s\bigr)(u,u)
    \geq
    sqK(-r)^s|u|^{2}\biggl(1-K\frac{Dns}{1-s}\biggr).
  \end{align*}
  Choosing $K=\frac{1-s}{2Dns}$ then proves the claim with
  \begin{align*}
  i\partial\dbar\bigl(-(-r)^{s}\bigr)(u,u)\geq \frac{(1-s)q}{4Dn}(-r)^{s}|u|^{2}.
  \end{align*}
\end{proof}

We are prepared to prove the main result of this section.
\begin{proposition}\label{P:dbarHardy}
Suppose $\Omega$ satisfies condition $(H_q)$. Let $s\in [0,1)$. 

Then there exists a smooth defining function, $r$, for $\Omega$ such that
\begin{equation}\label{E:dbarHardy}
  \int_{\Omega}(-r)^{-2+s}\left|\left[\dbarstar,r\right]u\right|^2
  \lesssim
  \bigl((-r)^{s}\dbar u,\dbar u\bigr)+\bigl((-r)^{s}\dbarstar u,\dbarstar u\bigr)
\end{equation}
holds for all $u\in\mathcal{D}^{0,q}(\Omega)$; the constant in $\lesssim$ depends on $s$.
\end{proposition}
\begin{proof}
  Let $s\in[0,1)$ be fixed and let $r$ be given by Lemma \ref{L:wlog}.
   Write $(-r)^{-2+s}=\frac{1}{1-s}\frac{\partial}{\partial r}(-r)^{-1+s}$ and express $\frac{\partial}{\partial
   r}$ as a linear combination of $\frac{\partial}{\partial\bar{z}_k}$ derivatives and vector fields which are 
   tangent to $b\Omega$.
   
  Then
  \begin{align*}
    \int_{\Omega}(-r)^{-2+s}&\bigl|\bigl[\dbarstar,r\bigr]u\bigr|^2\\
   &=\frac{1}{1-s}\int_{\Omega}\biggl(\sum_{k=1}a_{k}\frac{\partial}{\partial \bar{z}_{k}}(-r)^{-1+s}+D_t(-r)^{-1+s}\biggr)
   \bigl|\bigl[\dbarstar,r\bigr]u\bigr|^2,
  \end{align*}
  where $D_t(r)=0$ on $b\Omega$.
 Thus we can write 
  $D_t(r(z))=O(r(z))$, and it follows that
  \begin{align*}
  \frac{1}{1-s}\int_{\Omega}\bigl(D_t(-r)^{-1+s}\bigr)\bigl|\bigl[\dbarstar,r\bigr]u\bigr|^2
  &=\int_{\Omega}(-r)^{-2+s}D_t(r)\bigl|\bigl[\dbarstar,r\bigr]u\bigr|^2\\
  &\lesssim\int_{\Omega}(-r)^{-1+s}\bigl|\bigl[\dbarstar,r\bigr]u\bigr|^{2}\\
  &\lesssim \epsilon\bigl\|(-r)^{-1+\frac{s}{2}}\bigl[\dbarstar,r\bigr]u\bigr\|^2+\frac{1}{\epsilon}\bigl\|(-r)^{\frac{s}{2}}u\bigr\|^2,
  \end{align*}
  where the last line holds by (sc)-(lc) inequality and the fact that $[\dbarstar,r]$
  is in $L^\infty$. Furthermore, integration by parts gives
  \begin{align*}
    \int_{\Omega}\biggl(\sum_{k=1}a_{k}\frac{\partial}{\partial \bar{z}_{k}}(-r)^{-1+s}\biggr)
   \bigl|\bigl[\dbarstar,r\bigr]u\bigr|^2
   =\int_{\Omega}(-r)^{-1+s}\sum_{k=1}^{n}\frac{\partial}{\partial\bar{z}_k}
   \Bigl(a_k\bigl|\bigl[\dbarstar,r\bigr]u\bigr|^2\Bigr).
  \end{align*}
  Here the boundary integral vanishes since $[\dbarstar,r]u$ vanishes on the boundary. It follows that
  \begin{align*}
    \int_{\Omega}&(-r)^{-1+s}\sum_{k=1}^{n}\frac{\partial}{\partial\bar{z}_k}
    \Bigl(a_k\bigl|\bigl[\dbarstar,r\bigr]u\bigr|^2\Bigr)\\
   & \lesssim
    \int_{\Omega}(-r)^{-1+s}\bigl|\bigl[\dbarstar,r\bigr]u\bigr|^2
    +\int_{\Omega}(-r)^{-1+s}\Biggl|\biggl\langle   
    \sum_{k=1}^{n}\frac{\partial}{\partial\bar{z}_{k}}([\dbarstar,r]u),[\dbarstar,r]u\biggr\rangle\Biggr|\\
    &\lesssim \epsilon\int_{\Omega}(-r)^{-2+s}\bigl|\bigl[\dbarstar,r\bigr]u\bigr|^2+
    \frac{1}{\epsilon}\int_{\Omega}(-r)^{s}\biggl(|u|^2+\psum_{|I|=q}\sum_{k=1}^{n}\bigl|\frac{\partial  
    u_{I}}{\partial
    \bar{z}_{k}}\bigr|^2\biggr)
  \end{align*}
  Collecting the above estimates and choosing $\epsilon>0$ sufficiently small, we obtain
  \begin{align}\label{E:IntEst}
     \int_{\Omega}(-r)^{-2+s}\bigl|\bigl[\dbarstar,r\bigr]u\bigr|^2
     \lesssim
    \bigl\|(-r)^{\frac{s}{2}}u\bigr\|^2+\psum_{|I|=q}\sum_{k=1}^{n}\Bigl\|(-r)^{\frac{s}{2}}\frac{\partial u_{I}}
     {\partial\bar{z}_{k}}\Bigr\|^2. 
  \end{align}
  If $s=0$, then the left hand side of (\ref{E:IntEst}) is dominated by $\|\dbar u\|^2+\|\dbarstar u\|^2$, which  
  yields inequality (\ref{E:dbarHardy}).\\
  
   Now suppose $s\in(0,1)$. We recall that Lemma \ref{L:wlog} implies that 
  $(-r)^{s}|u|^2\lesssim i\partial\dbar(-(-r)^s)(u,u)$ holds. Thus, using Proposition \ref{P:BasEst} with 
  $\tau=(-r)^s$, we get
  \begin{align*}
     \int_{\Omega}(-r)^{-2+s}&\bigl|\bigl[\dbarstar,r\bigr]u\bigr|^2
     \lesssim -\int_{\Omega}i\partial\dbar(-r)^s(u,u)
     +\psum_{|I|=q}\sum_{k=1}^{n}\Bigl\|(-r)^{\frac{s}{2}}\frac{\partial u_{I}}{\partial\bar{z}_{k}}\Bigr\|^2\\
     &\lesssim
     \bigl\|(-r)^{\frac{s}{2}}\dbar u\bigr\|^2+\bigl\|(-r)^{\frac{s}{2}}\dbarstar u\bigr\|^2
     -2\text{Re}\bigl(\bigl[\dbarstar,(-r)^{s}\bigr]u,\dbarstar u\bigr).
  \end{align*}
  The last term on the right hand side can be easily controlled, in fact
  \begin{align*}
    \left|2\text{Re}\bigl(\bigl[\dbarstar, (-r)^s\bigr]u,\dbarstar u\bigr)\right|
    &\lesssim\bigl|\bigl((-r)^{\frac{s}{2}-1}[\dbarstar,r]u,(-r)^{\frac{s}{2}}\dbarstar u\bigr)\bigr|\\
    &\lesssim\epsilon\int_{\Omega}(-r)^{-2+s}\bigl|\bigl[\dbarstar,r\bigr]u\bigr|^2+\frac{1}{\epsilon}
   \bigl\|(-r)^{\frac{s}{2}}\dbarstar u\bigr\|^2,
  \end{align*}
  by the (sc)-(lc) inequality. Choosing $\epsilon>0$ sufficiently small, we obtain
  \begin{align*}
    \int_{\Omega}(-r)^{-2+s}\bigl|\bigl[\dbarstar,r\bigr]u\bigr|^2
     \lesssim
     \bigl\|(-r)^{\frac{s}{2}}\dbar u\bigr\|^2+\bigl\|(-r)^{\frac{s}{2}}\dbarstar u\bigr\|^2.
  \end{align*}
\end{proof}

\section{The Proof}\label{S:TheProof}

We state a quantitative form of Theorem \ref{T:MainT}:

\begin{theorem}\label{T:MainTheorem}
  Let $\Omega$ be a smoothly bounded domain satisfying condition $(H_q)$.
   Then the Bergman projection, $B_j$, is continuous on the Sobolev space $H_{0,j}^{s}(\Omega)$, $s>0$, for $j\in\{q-1,\dots,n-1\}$.
\end{theorem}
\begin{proof}
  We shall prove that $\|B_j f\|_{k}\leq C_{k}\|f\|_{k}$ holds for all integer $k>0$; the general case follows  
  by the usual interpolation arguments.
  
  The proof goes via a downward induction on $j$, the form level, as well as an upward induction on the 
  $k$, the order of differentiation, in the following manner.
  The induction basis (on the form level $j=n-1$) is satisfied: since the $\dbar$-Neumann problem  
  on  $(0,n)$-forms is an elliptic boundary value problem, the $\dbar$-Neumann operator, $N_n$, gains 
  two  derivatives, which implies that $\|B_{n-1}f\|_k\leq C_{k}\|f\|_k$ holds for all $k\geq0$ since
  $B_{n-1}=Id-\dbarstar N_n\dbar$.
  The induction basis (on the order of differentiation $k=0$) is also satisfied: $\|B_j f\|\leq 
  C_{0}\|f\|$ holds for all $j\in\{0,\dots,n-1\}$ by definition.
  
  In the following, we shall prove the case $k=1$ only, so that the main ideas are not cluttered by  
  technicalities. We indicate at the end how to prove the induction step for $k>1$.\\
  
  Let $j\in\{q,\dots,n-1\}$ be fixed. Suppose that $B_j$ is continuous on $H_{0,j}^{1}(\Omega)$. 
  We want to show that
  \begin{align}\label{E:1apriori}
    \|B_{j-1}f\|_1\leq C_{1}\|f\|_1
  \end{align} 
  holds for all $f\in \Lambda^{0,j-1}(\overline{\Omega})$. We first assume that $B_{j-1}f 
  \in\Lambda^{0,j-1}(\overline{\Omega})$ and establish \eqref{E:1apriori}. At the end of the proof, we show
  how to pass from this apriori estimate to a true estimate.\\
 
  \medskip
  
  \noindent{\bf I. Standard reduction.}
  \medskip 
  
  Let $r$ be a defining function for $\Omega$ which satisfies (\ref{E:dbarHardy}) for some
  $s=s_0\in(0,1)$ fixed, which will be chosen later.  We can 
  assume that 
  $\sum_{k=1}^{n}|r_{z_k}|^2\neq 0$ on a strip near $b\Omega$, i.e., on $S_{\eta}=\{z\in\Omega\;|\;-
  \eta<r(z)<0\}$ for some fixed $\eta>0$. Let $\chi$ be a smooth, non-negative function which vanishes  
  on
  $\Omega\backslash S_{2\eta}$ and equals $(\sum_{k=1}^{n}|r_{z_{k}}|^2)^{-1}$ on 
  $\bar{S}_{\eta}$. We define $(1,0)$-vector fields as follows:
  \begin{align*}
    L_j=\frac{\partial}{\partial z_j}-\chi\cdot r_{z_{j}}
    \sum_{l=1}^{n}r_{\bar{z}_{l}}\frac{\partial}{\partial z_l}\;\;
    \text{for}\;\;j\in\{1,\dots,n\}\;\;\text{and}\;\;
    \sn=\sum_{k=1}^{n}r_{\bar{z}_k}\frac{\partial}{\partial z_k}.
  \end{align*} 
   Then $\left\{L_1,\dots,L_{n-1}, \sn\right\}$ is a basis of $(1,0)$ vector fields. Also, define the $(1,1)$ vector field  
   $X=\sn-\bar\sn$. Notice that $L_j$, $1\leq j\leq n$, and $X$ are tangential. 
   Since the complex $\dbar\oplus\vartheta$ ($\vartheta$ being the formal adjoint of $\dbar$) is elliptic, the Sobolev $1$-norm of $B_{j-1}f$ is dominated 
   by  the following terms: the $L^2$-norm of $\dbar B_{j-1}f$ (which is $0$), the $L^2$-norm of $\vartheta 
   B_{j-1}f$ (which equals $\vartheta f$), the $L^2$-norm of $B_{j-1}f$, and the $L^2$-norms of tangential  
   derivatives of $B_{j-1}f$. That is
   \begin{align*}
     \|B_{j-1} f\|_{1}^2\lesssim\|f\|_{1}^{2}+\|B_{j-1}f\|^{2}+\sum_{k=1}^{n-1}\left(\|L_{k}B_{j-1}f\|^2
     +\|\bar{L}_{k}B_{j-1}f\|^2\right)+\|XB_{j-1}f\|^2.
   \end{align*}
   Let us see how to estimate terms involving barred derivatives of $B_{j-1}f$. First note that barred 
   derivatives of $B_{0}f$ vanish, since $B_{0}f$ is holomorphic. For $j\geq2$ we would like to use 
   inequality (\ref{E:VeryBasEst}), but $B_{j-1}f$ is not necessarily in $\mathcal{D}^{0,j-1}(\Omega)$.
   However, $f-B_{j-1}f$ is perpendicular to $\ker\dbar$, which implies that $f-B_{j-1}f$ is in
   $\mathcal{D}^{0,j-1}(\Omega)$ (see (\ref{D:DomdbarstarDef})). 
    Hence we can apply inequality (\ref{E:VeryBasEst}) to the term $f-B_{j-1}f$, which yields
    \begin{align}\label{E:barder}
      \|\bar{\sn}B_{j-1}f\|^{2}+\sum_{k=1}^{n}\|\bar{L}_{k}B_{j-1}f\|^{2}
      \lesssim
      \|f\|_{1}^{2}+\|\dbar(f-B_{j-1}f)\|^{2}
      \lesssim\|f\|_{1}^{2}.
    \end{align}
    To estimate the terms involving the tangential $(1,0)$ vector fields we use the standard argument that 
    integration by parts twice gives the following
    \begin{align*}
      \sum_{k=1}^{n}\|L_k B_{j-1}f\|^2
      \lesssim
      \sum_{k=1}^{n}([\bar{L}_k,L_k] B_{j-1}f, B_{j-1}f)+\sum_{k=1}^{n}\|\bar{L}_{k}B_{j-1}f\|^{2}
      +\|B_{j-1}f\|^2.
    \end{align*}
    Since the commutator $[\bar{L}_k,L_k]$, $1\leq k\leq n$, represents a tangential vector field, we can  
    express it as a linear combination of $X$, $L_l$ and $\bar{L}_l$ for $1\leq l\leq n$.  By inequality 
    (\ref{E:barder}) and the (sc)-(lc) inequality, it follows that
    \begin{align*}
       \sum_{k=1}^{n}\|L_k B_{j-1}f\|^2
       \lesssim
       \epsilon\|XB_{j-1}f\|^{2}+\frac{1}{\epsilon}\|B_{j-1}f\|^{2}+\|f\|_1^2
    \end{align*}
    holds for $\epsilon>0$.
    In particular, we obtain
    \begin{align}\label{E:OneNormB}
      \|B_{j-1}f\|_{1}^2\lesssim
      \|XB_{j-1}f\|^{2}+\|f\|_{1}^{2}+\|B_{j-1}f\|^{2}
      \lesssim
      \|XB_{j-1}f\|^{2}+\|f\|_{1}^{2},
    \end{align}
    where the last step results from the boundedness of $B_{j-1}$ on $L^{2}_{0,j-1}(\Omega)$. Thus it 
    remains to show that $\|XB_{j-1}f\|$ is dominated by the Sobolev $1$-norm of $f$.    
    
     \medskip
  
  \noindent{\bf II. Main apriori estimate (heart of the proof).}
  \medskip 
  
    From here on,  we 
    call those terms \emph{allowable} which are dominated by $\|f\|_{1}^{2}$. 
    In order to show that $\|XB_{j-1}f\|^{2}$ is allowable, we need to introduce a new vector field, $T$,  
    which equals $X$ at the first order level, but preserves membership in the domain of $\dbarstar$. 
    
    \begin{lemma}\label{L:T} There exists a smooth, tangential vector field $T$ such that
    \begin{itemize}
    \item[(i)] $T u\in\mathcal{D}^{0,j}(\Omega)$ 
    whenever $u\in\mathcal{D}^{0,j}(\Omega)$, and
    \item[(ii)] $\| (X-T)u\|\lesssim \|u\|$ for $\Lambda^{0,j}(\overline{\Omega})$.
    \end{itemize}
    \end{lemma}
    
    \begin{proof}
     Recall that we chose $\chi$ to be a smooth, non-negative function which vanishes on
    $\Omega\backslash S_{2\eta}$ and equals $(\sum_{k=1}^{n}|r_{z_{k}}|^2)^{-1}$ on $\bar{S}_{\eta}$ for 
    some fixed $\eta>0$.
     Set
    \begin{align}\label{D:tildeT}
      T u:=Xu+\chi\cdot\dbar r\wedge([X,[\dbarstar,r]]u).
    \end{align}
    
    Note that $T$ acts diagonally at the first order level, since $X$ does, and preserves the form level.
    It is straightforward to check that $T u\in\mathcal{D}^{0,j}(\Omega)$ 
    whenever $u\in\mathcal{D}^{0,j}(\Omega)$.
    The claimed property (ii), on $(X-T) u$, is obvious.
    \end{proof}
    
     Thus, to show that $\|XB_{j-1}f\|^{2}$ is 
    allowable, it suffices to prove the
    \begin{center}
      \underline{Claim}: $\|T B_{j-1}f\|^2$ is allowable.
    \end{center}
    We define $\varphi=N_{j}\dbar f$, note that then $\varphi\in\mathcal{D}^{0,j}(\Omega)$ and 
    $B_{j-1}f=f- 
    \dbarstar\varphi$. In order to deal with 
    certain error terms involving $\varphi$ arising in the proof of the Claim, we need the following lemma.
    \begin{lemma}\label{L:Estphi}
      \begin{itemize}
        \item[(i)] $Q(\varphi,\varphi)$ is allowable,
        \item[(ii)] $\|\varphi\|^{2}$, $\|\bar{\sn}\varphi\|^{2}$ and $\|\bar{L}_{k}\varphi\|^{2}$, $1\leq k\leq n$,  
        are allowable,
          \item[(iii)] $\|L_{k}\varphi\|^2\lesssim\epsilon\|T\varphi\|^{2}+\frac{1}{\epsilon}\|f\|^{2}$
           for $1\leq k\leq n$, $\epsilon>0$,
        \item[(iv)] $\|T\varphi\|^{2}\lesssim\|f\|_{1}^{2}+\|B_{j-1}f\|_{1}^{2}$.
      \end{itemize}
    \end{lemma}
    \begin{proof}[Proof of Lemma \ref{L:Estphi}]
    The proof of (i) follows directly from the definition of $\varphi$, i.e.,
    \begin{align*}
      Q(\varphi,\varphi)&=(\dbar\varphi,\dbar\varphi)+(\dbarstar\varphi,\dbarstar\varphi)\\
      &=(\dbar N_j\dbar f, \dbar N_j \dbar f)+(\dbarstar N_j\dbar f,\dbarstar N_j\dbar f)\\
      &=\|f-B_{j-1} f\|^{2}\lesssim\|f\|^{2}.
    \end{align*}
    
    To prove (ii), we use inequality (\ref{E:VeryBasEst}) and (i), that is
    \begin{align*}
      \|\varphi\|^{2}+\|\bar{L}_{k}\varphi\|^{2}\lesssim Q(\varphi,\varphi)\lesssim\|f\|^2.
    \end{align*}
    
    Inequality (iii) follows from the same arguments as those made directly below inequality 
    (\ref{E:barder}).
    
    For the proof of (iv) we use the Boas-Straube formula in \cite{Boa-Str90}, which expresses 
    $\dbarstar N_j$ in terms of $B_{j-1}$, $B_j$ and $N_{t,j}$; here $N_{t,j}$ is the solution operator to 
    the  weighted $\dbar$-Neumann problem with weight $w_t(z)=\exp(-t|z|^2)$. Recall that 
    $\varphi=N_j\dbar f$, and hence we are interested in the operator $N_j\dbar$ and not in the operator
     $\dbarstar N_j$. However, $\dbarstar N_j$ is the $L^2$-adjoint of $N_j\dbar$, and thus the formula 
     for $\dbarstar N_j$ in \cite{Boa-Str90}, pg. 29, implies
    \begin{align*}
      N_j\dbar f=B_j\{w_t N_{t,j}\dbar w_{-t}(f-B_{j-1}f)\}.
    \end{align*} 
    The induction hypothesis says that $B_j$ is continuous on $H_{0,j}^{1}(\Omega)$. Thus
    \begin{align*}
      \|T\varphi\|^2=\|T N_{j}\dbar f\|^2\lesssim\|N_{t,j}\dbar (w_{-t}f-w_{-t}B_{j-1}f)\|_{1}^2.
    \end{align*}
    
    A theorem of Kohn in \cite{Koh73} 
    implies that $N_{t,j}\dbar$ is continuous on 
    $H_{0,j-1}^{1}(\Omega)$ as long as $t>0$ is sufficiently large.
    Actually, in \cite{Koh73}, Kohn assumes that $b\Omega$ is
    pseudoconvex, so his result is not immediately applicable under our hypotheses.  
    However, on $(0,j)$-forms the basic $\bar\partial$-Neumann identity is (\ref{E:BasId}) and
    the boundary integrand is non-negative by condition ($H_j$). The techniques in  \cite{Koh73}
    may now be applied to give the claimed estimate on $N_{t,j}\dbar$.
    
    Thus we obtain
    \begin{align*}
      \|T\varphi\|^{2}\lesssim\|f\|_{1}^{2}+\|B_{j-1}f\|_{1}^{2},
    \end{align*}
    which proves (iv).
    \end{proof}
    
    Now we are ready to show that $\|T B_{j-1}f\|^{2}$ is allowable. Since $B_{j-1}f=f-\dbarstar\varphi$, it  
    follows 
    that
    \begin{align*}
     \left\|T B_{j-1}f\right\|^{2}
      &=\left(T f,T B_{j-1}f\right)-\left(T\dbarstar\varphi,T B_{j-1}f\right)\\
      &=\left(T f,T B_{j-1}f\right)-\left(\varphi,\left[\dbar,T^*T\right] B_{j-1}f\right)-
      \left(\varphi,
      T^*T\dbar B_{j-1}f\right).
    \end{align*}
    Since $\dbar B_{j-1}f$ equals $0$, it follows with the (sc)-(lc) inequality
    \begin{align*}
      \left\|T B_{j-1}f\right\|^{2}\lesssim\left\|f\right\|_{1}^{2}
      +\left|\left(\left[\dbar,T^*T\right]\varphi,B_{j-1}f
      \right)\right|.
    \end{align*}
    Also, 
    \begin{align*}
      \left[\dbar,T^{*}T\right]=\left[\dbar,T^*\right] T+\left[\dbar,T\right]T^*
      -\left[\left[\dbar,T\right],T^*\right],
    \end{align*}
    where the last term above is of order $1$.  The adjoint of $[\dbar,T^*]$  is $[T,\dbarstar]$, 
    but the adjoint of $[\dbar,T]$ is not $[T^{*},\dbarstar]$, since $T^*$ does not preserve the domain of 
    $\dbarstar$. However, since $X$ is self-adjoint, $T$ and $T^*$ only differ by terms of order $0$, which 
    allows us to integrate by parts with negligible error terms. That is,
    \begin{align*}
      \left(\varphi,\left[\dbar,T\right]T^*B_{j-1}f\right)=\left(T^*\dbarstar\varphi,T^*B_{j-1}f\right)
      -\left(T^*\varphi,\dbar T^*B_{j-1}f\right),
    \end{align*}
    and
    \begin{align*}
      \left(T^*\varphi,\dbar T^*B_{j-1}f\right)
      &\leq
      \left(T\varphi,\dbar T^*B_{j-1}f\right)+\left\|\varphi\right\|\cdot\left\|\dbar T^*B_{j-1}f\right\|\\
      &=\left(\dbarstar T\varphi,T^*B_{j-1}f\right)+\left\|\varphi\right\|\cdot
      \left\|\left[\dbar,T^*\right]B_{j-1}f\right\|.
    \end{align*}
   Thus, using again that $T$ and $T^*$ are equal at the first order level, we get
    \begin{align*}
      \left\|T B_{j-1}f\right\|^{2}\lesssim
      \left|\left(\left[T,\dbarstar\right]\varphi,TB_{j-1}f\right)\right|
      +\left\|f\right\|_{1}^{2}+\left\|\varphi
      \right\|\cdot\left\|B_{j-1}f\right\|_{1},
    \end{align*}
    where the last step follows since $B_{j-1}f$ is in the kernel of $\dbar$.
    Inequality (\ref{E:OneNormB}), part (ii) of Lemma \ref{L:Estphi} and the (sc)-(lc) inequality yield
    \begin{align}\label{E:toshowT}
      \left\|T B_{j-1}f\right\|^{2}\lesssim\left\|\left[T,\dbarstar\right]\varphi\right\|^{2}+\left\|f\right\|_{1}^{2}.
    \end{align}
    Therefore, our claim follows by the (sc)-(lc) inequality  if $\left\|\left[T,\dbarstar\right]\varphi\right\|^2$ is  
    allowable. 
    
    The term $[T,\dbarstar]\varphi$ is equal to $[X,\vartheta]\varphi$ up to terms of order $0$, where 
    $\vartheta$ is the formal adjoint of $\dbar$. A straighforward computation gives
    \begin{align*}
      \left[X,\vartheta\right]\varphi=\psum_{|J|=j-1}\sum_{k,l=1}^{n}\left\{\frac{\partial^2 r}{\partial  
      z_{l}\bar{z}_{k}}
      \frac{\partial\varphi_{lJ}}{\partial z_{k}}
      -\frac{\partial^{2} r}{\partial z_{l}\partial z_{k}}\frac{\partial\varphi_{lJ}}{\partial\bar{z}_{k}}\right\}
      d\bar{z}^{J}.
    \end{align*}
    Note that the terms involving $\frac{\partial\phi_{lJ}}{\partial\bar z_k}$ are allowable by 
    part (ii) of Lemma \ref{L:Estphi}. To deal with the remaining terms, notice that
    \begin{align*}
      \frac{\partial}{\partial z_{k}}=L_{k}+\chi r_{z_{k}}X+ 
      \chi r_{z_{k}}\sum_{i=1}^{n}r_{z_{i}}\frac{\partial}{\partial\bar{z}_{k}}\;\;\text{for}\;\;k\in\{1,\dots,n\}.   
    \end{align*}
    Thus we obtain
    \begin{align*}
       \psum_{|J|=j-1}\sum_{k,l=1}^{n}\frac{\partial^2 r}{\partial z_{l}\bar{z}_{k}}
      \frac{\partial\varphi_{lJ}}{\partial z_{k}}
     =\psum_{|J|=j-1}\sum_{k,l=1}^{n}\left\{\frac{\partial^{2}r}{\partial z_{l}\partial\bar{z}_{k}}\right.
     ((\chi X\varphi_{lJ})r_{z_{k}}&+L_{k}\varphi_{lJ})
     \left. +\chi r_{z_{k}}\sum_{i=1}^{n}r_{z_{i}}\frac{\partial\varphi_{lJ}}{\partial \bar{z}_{k}}\right\}.
    \end{align*}
    By part (ii) and (iii) of Lemma \ref{L:Estphi},  applied to the last two terms of the right hand side of the 
    above equation, it follows that
    \begin{align}\label{E:beforeCSprep}
      \left\|\left[T,\dbarstar\right]\varphi\right\|^{2}
      \lesssim
      \left\|\psum_{|J|=j-1}\sum_{k,l=1}^{n}\frac{\partial^{2}r}{\partial z_{l}\partial\bar{z}_{k}}
     (T\varphi_{lJ})r_{z_{k}}\right\|^{2}+\text{ allowable}.
    \end{align}
    \begin{lemma}\label{L:CSprep}
      There exists a $(0,j)$-form $\psi$ such that
      \begin{align}\label{E:CSprep}
        \left|\psum_{|J|=j-1}\sum_{k,l=1}^{n}\frac{\partial^{2}r}{\partial z_{l}\partial\bar{z}_{k}}
     (T\varphi_{lJ})r_{z_{k}}\right|
     \lesssim |i\partial\dbar r(T\varphi,\psi)|+\left|\left[\dbarstar,r\right]T\varphi\right|.
      \end{align}
    \end{lemma}
    \begin{proof}
      For notational ease, let us write $\phi$ for $T\varphi$ temporarily. We define 
      \begin{align*}
        \psi=\dbar\biggl(r\psum_{|J|=j-1}d\bar{z}^{J}\biggr).
      \end{align*}  
       For $j=1$,  (\ref{E:CSprep}) holds trivially, since $\overline{\psi}_{k}=r_{z_{k}}$. 
       For $j>1$, an error term occurs when passing to $\psi$. 
       In the following, we shall indicate what kind of
       algebraic manipulations of this error term  lead to \eqref{E:CSprep}.  For that we need to fix some 
       notation:  Recall that we write $kI$ for $\{k,I\}$, if $I$ is an increasing multi-index  and  
       $k\notin I$. Furthermore, we shall mean by $k\cup I$  the increasing multi-index which equals $kI$ 
       as a set. As before, $\sigma^{I}_{J}$ is the sign of the permutation 
       $\left(
       \begin{smallmatrix}
       I\\ J
       \end{smallmatrix}
       \right)$
       and is zero whenever $I$ 
       and $J$ are not equal as set.
       
       Notice first 
       that 
       $ \overline{\psi}_{I}=\sum_{m\in I}r_{z_{m}}\sigma^{m(I\backslash m)}_{I}$
       for any increasing multi-index $I$ of length $j$.
       Moreover, if $J$ is an increasing multi-index of length $j-1$ and $k\notin J$, we can write
       \begin{align*}
         \overline{\psi}_{kJ}=r_{z_{k}}+\sum_{m\in J}r_{z_{m}}\sigma^{m(k\cup(J\backslash m))}_{kJ}.
       \end{align*}
       Using this, a straightforward computation then gives
       \begin{align*}
         \psum_{|J|=j-1}\sum_{k,l=1}^{n}&\frac{\partial^2 r}{\partial z_{l}\partial\bar{z}_{k}}
         \phi_{lJ}r_{z_{k}}-i\partial\dbar r(\phi,\psi)\\
         =&\sum_{k,l=1}^{n}\frac{\partial^2 r}{\partial 
         z_{l}\partial\bar{z}_{k}}
         \left\{\mathop{\psum_{|J|=j-1}}_{k\in J}\phi_{lJ}r_{z_{k}}
        -\mathop{\psum_{|J|=j-1}}_{k\notin J}
         \phi_{lJ}\Biggl[\sum_{m\in J}
         r_{z_{m}}\sigma^{m(k\cup(J\backslash m))}_{kJ}\Biggr]\right\}.
       \end{align*}   
       Let us consider the two terms in the parentheses on the right hand side for  $k=l$ fixed. Notice that in 
       this case the first term vanishes, so we only need to study the second term. That 
       is
       \begin{align*}
         \psum_{|J|=j-1}
         \phi_{kJ}\Biggl[\sum_{m\in J}
         r_{z_{m}}\sigma^{m(k\cup(J\backslash m))}_{kJ}\Biggr]
         &=\psum_{|I|=j-2}\sum_{m=1}^{n}
      r_{z_{m}}\phi_{k(m\cup I)}\cdot\sigma^{k(m\cup I)}_{m(k\cup I)}\\
      &=\psum_{|I|=j-2}\sum_{m=1}^{n}
      r_{z_{m}}\phi_{m(k\cup I)}
      =\mathop{\psum_{|J|=j-1}}_{k\in J}\sum_{m=1}^{n}r_{z_{m}}\phi_{mJ},
       \end{align*}
       which equals the sum over those components of $-[\dbarstar,r]\phi$ whose multi-indices contain $k$.
       With similar, though more elaborate, computations one obtains for the case where $k\neq l$ are 
       fixed the following
       \begin{align*}
         \mathop{\psum_{|J|=j-1}}_{k\in J}\phi_{lJ}r_{z_{k}}
        -\mathop{\psum_{|J|=j-1}}_{k\notin J}
         \phi_{lJ}\Biggl[\sum_{m\in J}
         r_{z_{m}}\sigma^{m(k\cup(J\backslash m))}_{kJ}\Biggr]
         &=\mathop{\psum_{|J|=j-1}}_{l\in J,k\notin J}\sigma^{l(k\cup(J\backslash l))}_{kJ}
         \sum_{m=1}^{n}r_{z_{m}}\phi_{mJ}\\
         &=-\mathop{\psum_{|J|=j-1}}_{l\in J,k\notin J}\sigma^{l(k\cup(J\backslash l))}_{kJ}
         \left([\dbarstar,r]\phi\right)_{J}.
       \end{align*}
       Thus the error term appearing when passing to $\psi$ is described by the following equation
         \begin{align*}
            \psum_{|J|=j-1}\sum_{k,l=1}^{n}\frac{\partial^2 r}{\partial z_{l}\partial\bar{z}_{k}}\phi_{lJ}r_{z_{k}}
            -&
            i\partial\dbar r(\phi,\psi)\\
            =\sum_{k=1}^{n}\frac{\partial^2 r}{\partial 
            z_{k}\partial\bar{z}_{k}}
            \Biggl\{
            \mathop{\psum_{|J|=j-1}}_{k\in J}\left([\dbarstar,r]\phi\right)_{J}\Biggr\}
           &-
           \mathop{\sum_{k,l=1}^{n}}_{k\neq l}\frac{\partial^2 r}{\partial z_{l}\bar{z}_{k}} 
           \Biggl\{\mathop{\psum_{|J|=j-1}}_{l\in J,k\notin J}\sigma^{l(k\cup(J\backslash l))}_{kJ}
           \left([\dbarstar,r]\phi\right)_{J}
           \Biggr\},
         \end{align*}
         which implies the claimed estimate \eqref{E:CSprep}.  
    \end{proof}
    Recall that $i\partial\dbar r(u,u)\geq 0$ holds for all $(0,j)$-forms $u$. Hence, it follows by the 
    Cauchy-Schwarz inequality that
    \begin{align}\label{E:CauSch}
      |i\partial\dbar r(T\varphi,\psi)|\leq\left(i\partial\dbar r(T\varphi,T\varphi)\right)^{\frac{1}{2}}\cdot
      \left(i\partial\dbar r(\psi,\psi)\right)^{\frac{1}{2}}.
    \end{align}
    This, combined with Lemma \ref{L:CSprep} and \eqref{E:beforeCSprep}, implies that
    \begin{align}\label{E:MainInt}
      \left\|\left[T,\dbarstar\right]\varphi\right\|^{2}
      \lesssim\int_{\Omega}
      i\partial\dbar r(T\varphi,T\varphi)+\left\|\left[\dbarstar,r\right]T\varphi\right\|^{2}+\text{ allowable}.
    \end{align}
    Recall that $T u\in\mathcal{D}^{0,j}(\Omega)$   
     whenever $u\in\mathcal{D}^{0,j}(\Omega)$.
     Thus, we can apply \eqref{E:BasEst} to \eqref{E:MainInt}, with $\tau=(-r)$, 
     and obtain
    \begin{align*}
    \int_{\Omega}
      i\partial\dbar r(T\varphi,T\varphi)
      \lesssim
      &\left\|(-r)^{\frac{1}{2}}\dbar T\varphi\right\|^{2}+\left\|(-r)^{\frac{1}{2}}\dbarstar T\varphi\right\|^{2} \\
      &+\left|\left(\left[\dbarstar,r\right]T\varphi,\dbarstar T\varphi\right)\right|+\text{ allowable}.
    \end{align*}
    Using the (sc)-(lc) inequality, we estimate the last term on the right hand side above 
     as follows:
    \begin{align*}
      \left|\left(\left[\dbarstar,r\right]T\varphi,\dbarstar T\varphi\right)\right|
      \lesssim
      \int_{\Omega}(-r)^{-2+\frac{1}{2}}\left|\left[\dbarstar,r\right]T\varphi\right|^{2}
      +\left\|(-r)^{\frac{3}{4}}\dbarstar T\varphi\right\|^{2}.
    \end{align*}
    Since $r$ is bounded from below, this implies that
    \begin{align*}
      \left\|\left[T,\dbarstar\right]\varphi\right\|^{2}
      \lesssim
      \int_{\Omega}(-r)^{-2+\frac{1}{2}}\left|\left[\dbarstar,r\right]T\varphi\right|^{2}
     +\left\|(-r)^{\frac{1}{2}}\dbar T\varphi\right\|^{2}&+\left\|(-r)^{\frac{1}{2}}\dbarstar T\varphi\right\|^{2}\\
     &+\text{ allowable}.
    \end{align*}
    Recall that we chose $r$ to be a defining function of $\Omega$ which satisfies (\ref{E:dbarHardy}) for 
    some fixed $s_{0}\in(0,1)$. We now choose $s_{0}=\frac{1}{2}$, so that (\ref{E:dbarHardy}) gives
    \begin{align*}
      \int_{\Omega}(-r)^{-2+\frac{1}{2}}\left|\left[\dbarstar,r\right]T\varphi\right|^{2}
      \lesssim
     \left\|(-r)^{\frac{1}{4}}\dbarstar T\varphi\right\|^{2}+\left\|(-r)^{\frac{1}{4}}\dbar T\varphi\right\|^{2}.
    \end{align*}
     Again, since $r$ is
     bounded from below, we may more simply write
    \begin{align*}
      \left\|\left[T,\dbarstar\right]\varphi\right\|^{2}\lesssim\left\|(-r)^{\frac{1}{4}}\dbar T\varphi\right\|^{2}+
      \left\|(-r)^{\frac{1}{4}}
       \dbarstar T\varphi\right\|^{2}+\text{ allowable}.
    \end{align*}
    
    The two terms on the right hand side of the above inequality contain a factor $(-r)^{\frac{1}{4}}$, which  makes them possible to estimate. 
    Consider $\|(-r)^{\frac{1}{4}}DT\varphi\|^{2}$, where $D$ will 
    either be $\dbar$ or $\dbarstar$. Let $\epsilon>0$ be a fixed number, which will be chosen later. 
    Note that $(-r)^{\frac{1}{2}}<\epsilon$ holds on the set
    $S_{\epsilon^{2}}=\{z\in\Omega\;|\;-\epsilon^2<r(z)<0\}$. Let $\zeta$ be a smooth, non-negative function 
    such 
    that $\zeta=1$ on $\Omega\backslash S_{\epsilon^{2}}$ and $\zeta=0$ on 
    $S_{\frac{\epsilon^{2}}{2}}$. Then
    \begin{align*}
      \|(-r)^{\frac{1}{4}}DT\varphi\|^{2}\lesssim\|\zeta DT\varphi\|^{2}+
      \epsilon\|(1-\zeta)DT\varphi\|^{2}=:A_{1}+A_{2}.
    \end{align*}
    To show that the term $A_1$ is allowable we commute:
    \begin{align*}
      A_1\lesssim\left\|\zeta\left[D,T\right]\varphi\right\|^{2}+\left\|\zeta T D\varphi\right\|^{2}
      \lesssim
      \left\|\left[\zeta,[D,T]\right]\varphi\right\|^{2}+\left\|\zeta\varphi\right\|_{1}^{2}+
      \left\|\zeta T D\varphi\right\|^{2}.
    \end{align*}
    Since $\zeta\varphi$ is compactly supported in $\Omega$, interior elliptic estimates give us
    \begin{align*}
      \|\zeta\varphi\|_{1}^{2}\lesssim Q(\zeta\varphi,\zeta\varphi)\lesssim 
      \epsilon^{-4}Q(\varphi,\varphi)\lesssim\epsilon^{-4}\|f\|^2,
    \end{align*}
    where the last step follows by part (i) of Lemma \ref{L:Estphi}. Also, $\left[\zeta,[D,T]\right]$ is of order $0$. 
    Thus, with part (ii) of Lemma \ref{L:Estphi}, we obtain
    \begin{align*}
      A_{1}\lesssim \epsilon^{-4}\|f\|^{2}+\left\|\zeta T D\varphi\right\|^{2}.
    \end{align*}
    If $D=\dbar$, then the second term on the right hand side of the last inequality equals $0$. If 
    $D=\dbarstar$, then
    \begin{align*}
     \left \|\zeta T\dbarstar\varphi\right\|^{2}
      &\lesssim\left\|\zeta T B_{j-1}f\right\|^{2}+\left\|\zeta T f\right\|^{2}\\
      &\lesssim \left\|\left[\zeta,T\right]B_{j-1}f\right\|^{2}+\left\|\zeta B_{j-1}f\right\|_{1}^{2}+\|f\|_{1}^{2}.
    \end{align*}
    Again, $\zeta B_{j-1}f$ is supported in $\Omega$, thus
    \begin{align*}
      \left\|\zeta B_{j-1}f\right\|^{2}\lesssim Q\left(\zeta B_{j-1}f,\zeta B_{j-1}f\right)\lesssim \epsilon^{-4}\|f\|_{1}^{2}.
    \end{align*}
    This concludes the proof of the term $A_1$ being allowable.
    To estimate the term $A_2$ we commute 
    again:
    \begin{align*}
      A_{2}\lesssim\epsilon\left\|DT\varphi\right\|^{2}
      \lesssim
      \epsilon\left(\left\|T D\varphi\right\|^{2}+\left\|\left[D,T\right]\varphi\right\|^{2}\right)
      \lesssim
      \epsilon\left(\left\|T D\varphi\right\|^{2}+\left\|\varphi\right\|_{1}^{2}\right).
    \end{align*}
    The term $T D\varphi$ vanishes if $D=\dbar$, otherwise
    \begin{align*}
      \left\|T D\varphi\right\|^{2}=\left\|T\dbarstar\varphi\right\|^{2}
      \lesssim
      \left\|T B_{j-1}f\right\|^{2}+\left\|T f\right\|^{2}
      \lesssim
      \left\|B_{j-1}f\right\|_{1}^{2}+\left\|f\right\|_{1}^{2}.
    \end{align*}
    Therefore we get
    \begin{align*}
      A_{1}+A_{2}\lesssim \epsilon^{-4}\|f\|_{1}^{2}+\epsilon\left(\left\|\varphi\right\|_{1}^{2}+\left\|B_{j-1}f\right\|_{1}^{2}\right)
      \lesssim\epsilon^{-4}\|f\|_{1}^{2}+\epsilon\|B_{j-1}f\|_{1}^{2},
    \end{align*}
    here the last estimate holds by part (iv) of Lemma \ref{L:Estphi}. This implies in particular 
    \begin{align*}
      \left\|\left[T,\dbarstar\right]\varphi\right\|^{2}\lesssim\epsilon^{-4}\|f\|_{1}^{2}+\epsilon
      \left\|B_{j-1}f\right\|_{1}^{2}
      \lesssim \epsilon^{-4}\|f\|_{1}^{2}+\epsilon\left\|T B_{j-1}f\right\|^{2},
    \end{align*}
    where the last line follows by inequality (\ref{E:OneNormB}). Invoking inequality (\ref{E:toshowT}) we 
    at last obtain
    \begin{align*}
      \left\|T B_{j-1}f\right\|^{2}\lesssim\epsilon^{-4}\|f\|_{1}^{2}+\epsilon\left\|T B_{j-1}f\right\|^{2}.
    \end{align*}
    Choosing $\epsilon>0$ sufficiently small so that we can absorb the last term into the left hand side, we 
    obtain that $\|T B_{j-1}f\|^{2}$ is allowable.  Thus we have shown that \eqref{E:1apriori} holds,
    assuming both $f$ and $B_{j-1}f$ are in $\Lambda^{0,j-1}(\overline{\Omega})$. 
    
    \medskip
  
  \noindent{\bf III. Removing smoothness assumptions.}
  \medskip 
    Consider a sequence of approximating subdomains: 
    define $\Omega_{\delta}=\{z\in\mathbb{C}^{n}\;|\;\rho_{\delta}(z)=\rho(z)+\delta|z|^{2}<0\}$ for 
    sufficiently small $\delta>0$. Note that $\Omega_{\delta}$ is smoothly bounded and strongly 
    pseudoconvex. This implies in particular that the Bergman projection, $B_{j-1}^{\delta}$, on 
    $\Omega_{\delta}$ preserves $\Lambda^{0,j-1}(\overline{\Omega}_{\delta})$.
    We can execute the above argument, replacing $\rho$ by 
    $\rho_{\delta}$, for $B_{j-1}^{\delta}$ on $\Omega_{\delta}$ to obtain
    \begin{align*}
     \left\|B_{j-1}^{\delta}f\right\|_{1,\Omega_{\delta}}\leq 
     C_{1}\left\|f\right\|_{1,\Omega_{\delta}}\;\;\text{for}\;\;
     f\in\Lambda^{0,j-1}(\overline{\Omega}_{\delta}),
    \end{align*}
    where the constant $C_{1}$ does not depend on $\delta$ since none of our estimates do. Thus
    $\|B_{j-1}^{\delta}f\|_{1,\Omega_{\delta}}\leq C_{1}\|f\|_{1,\Omega}$ for 
    $f\in\Lambda^{0,j-1}(\overline{\Omega})$ is a true estimate, which holds uniformly in $\delta$. 
    Since $B_{j-1}^{\delta}f$ converges to $B_{j-1}f$ pointwise, it follows that $B_{j-1}$ is in fact in
    $H_{0,j-1}^{1}(\Omega)$ and $\|B_{j-1}f\|_{1,\Omega}\leq C_{1}\|f\|_{1,\Omega}$ for 
    $f\in\Lambda^{0,j-1}(\overline{\Omega})$. Because $\Lambda^{0,j-1}(\overline{\Omega})$ is dense in 
    $H_{0,j-1}^{1}(\Omega)$ with respect to the Sobolev $1$-norm, continuity of $B_{j-1}$ on 
    $H_{0,j-1}^{1}(\Omega)$ follows. 
    
    \medskip
  
  \noindent{\bf IV. The induction step for $k>1$. }
  \medskip 
  
     Let $j\in\{q,\dots,n-1\}$ be fixed. Suppose $B_{j}$ is continuous on 
    $H_{0,j}^{k}(\Omega)$, and $B_{j-1}$ is continuous on $H_{0,j-1}^{k-1}(\Omega)$. In the following, we 
    indicate how to obtain that $\|B_{j-1}f\|_{k}\leq C_{k}\|f\|_{k}$ holds for 
    $f\in\Lambda^{0,j-1}(\overline{\Omega})$ assuming that $B_{j-1}f\in\Lambda^{0,j
    -1}(\overline{\Omega})$. 
    With arguments similar to the ones 
    preceding inequality (\ref{E:OneNormB}), we obtain
    \begin{align*}
    \left\|B_{j-1}f\right\|_{k}^{2}\lesssim\left\|X^{k}B_{j-1}f\right\|^{2}+\|f\|_{k}^{2}+\left\|B_{j-1}f\right\|_{k
    -1}^{2}
    \lesssim\left\|X^{k}B_{j-1}f\right\|^{2}+\|f\|_{k}^{2}\tag{$\widetilde{\ref{E:OneNormB}}$},
    \end{align*}
    where the last step follows by the induction hypothesis.
    For more details on how to derive ($\widetilde{\ref{E:OneNormB}}$) see also \cite{Boa-Str91}, pg. 83. 
    
    Thus one needs to show that $\|X^{k}B_{j-1}f\|$ is dominated by the Sobolev $k$-norm of $f$. In this 
    context 
    we call a term allowable if it is dominated by $\|f\|_{k}^{2}$. 
    As before, it is sufficient to show that $\|T^k B_{j-1}f\|^{2}$ is allowable.
    Again, we define $\varphi=N_{j}\dbar f$, 
    and compute
    \begin{align*}
      \left\|T^{k} B_{j-1}f\right\|^{2}=&\left(T^kf,T^k B_{j-1}f\right)-\left(T^k\dbarstar\varphi,T^k B_{j-1}f\right)\\
      =&\left(T^k f, T^k B_{j-1}f\right)-\left(\varphi,\left[\dbar,(T^{k})^*T^{k}\right]B_{j-1}f\right).
    \end{align*}
    Thus (sc)-(lc) inequality gives
    \begin{align*}
      \left\|T^k B_{j-1} f\right\|^{2}
      \lesssim
      \left\|f\right\|_{k}^{2}+\left|(\varphi,\left[\dbar,(T^{k})^{*}T^k\right]B_{j-1}f\right).
    \end{align*}
    Using that $(T^k)^*=(T^*)^{k}$ a straightforward computation gives
    \begin{align*}
      \left[\dbar,(T^k)^{*}T^k\right]=&
      3(T^*)^{k-1}\left[\dbar,T\right]T^{k}+3T^{k-1}\left[\dbar,T^{*}\right](T^*)^{k}\\&+
      \text{terms of order }2k-1.   
    \end{align*}
    Therefore
   \begin{align*}
      \left\|T^{k}B_{j-1}f\right\|^2\lesssim
      \left|\left(\left[T,\dbarstar\right]T^{k-1}\varphi,T^{k}B_{j-1}f\right)\right|
      +\left\|f\right\|_{k}^{2}+\left\|\varphi\right\|_{k-1}\left\|B_{j-1}f\right\|_{k}.
   \end{align*}    
   Since $B_{j-1}f$ is continuous on $H_{0,j-1}^{k-1}(\Omega)$, it follows that $\|\varphi\|_{k-1}^2$ is  
   allowable.
   Thus, using the (sc)-(lc) inequality and ($\widetilde{\ref{E:OneNormB}}$), it follows
   \begin{align*}
     \left\|T^k B_{j-1}f\right\|^{2}\lesssim
     \left\|\left[T,\dbarstar\right]T^{k-1}\varphi\right\|^{2}
     +\left\|f\right\|_{k}^{2}.\tag{$\widetilde{\ref{E:toshowT}}$}
   \end{align*}
  Obviously, one now needs estimates  
   for $\varphi$ similar to those estimates in Lemma \ref{L:Estphi}. By 
   analogous arguments as those which give Lemma \ref{L:Estphi}, using the induction hypotheses, one 
   obtains
   \begin{proof}[\textup{\textbf{Lemma $\widetilde{\textbf{\ref{L:Estphi}}}$}}]
     \begin{itemize}
       \item[($\tilde{\text{i}}$)]\emph{$Q(T^{k-1}\varphi,T^{k-1}\varphi)$ is allowable,}
       \item[($\widetilde{\text{ii}}$)]\emph{$\|\varphi\|_{k-1}^{2}$, $\|\bar{\sn}\varphi\|^2_{k-1}$ 
       and $\|\bar{L}_{l}\varphi\|_{k-1}^{2}$, 
       $1\leq l\leq n$, are allowable,}
       \item[($\widetilde{\text{iii}}$)]\emph{$\|L_{l}\varphi\|_{k-1}^{2}\lesssim\epsilon\|T\varphi\|_{k 
       -1}^{2}+\frac{1}{\epsilon}\|f\|_{k}^{2}$ for $1\leq l\leq n$, $\epsilon>0$,}
       \item[($\widetilde{\text{iv}}$)]\emph{$\|T\varphi\|_{k-1}^{2}\lesssim\|f\|_{k}^{2}+\|B_{j-1}f\|_{k}^{2}$.}
     \end{itemize}
     \renewcommand{\qedsymbol}{}
   \end{proof}
   Now one just follows the proof for $k=1$ starting at (\ref{E:toshowT}), with $T^{k-1}\varphi$ in place of 
   $\varphi$ and using Lemma $\widetilde{\ref{L:Estphi}}$ instead of Lemma \ref{L:Estphi}. This leads to
   the estimate 
   \begin{align*}
    \left\|T^k B_{j-1}f\right\|^2\lesssim\left\|f\right\|_{k}^{2},
   \end{align*}
   from which the apriori estimate $\|B_{j-1}f\|_{k}\leq C_{k}\|f\|_{k}$ follows, under the
   assumptions that $f, B_{j-1}f\in\Lambda^{0,j-1}(\overline{\Omega})$.  Passing from these
   apriori estimates to actual estimates follows as before.
 \end{proof}

\end{document}